\documentclass[10pt,letterpaper,dvips]{article}

 \usepackage{amssymb,latexsym,amsmath,amsthm}

\usepackage{fullpage}

\newtheorem*{thm*}{Theorem A}

\newtheorem{thm}{Theorem}

\newtheorem{lemma}{Lemma}

\author{A. Aghajani\thanks{School of Mathematics, Iran University of Science and Technology, Narmak, Tehran, Iran. Email: aghajani@iust.ac.ir.} \and
    C. Cowan\thanks{Department of Mathematics, University of Manitoba, Winnipeg, Manitoba, Canada R3T 2N2. Email: craig.cowan@umanitoba.ca. Research supported in part by NSERC.} }

\begin{document}

\def\d{ \partial_{x_j} }
\def\Na{{\mathbb{N}}}

\def\Z{{\mathbb{Z}}}

\def\IR{{\mathbb{R}}}

\newcommand{\E}[0]{ \varepsilon}

\newcommand{\la}[0]{ \lambda}

\newcommand{\s}[0]{ \mathcal{S}}

\newcommand{\AO}[1]{\| #1 \| }

\newcommand{\BO}[2]{ \left( #1 , #2 \right) }

\newcommand{\CO}[2]{ \left\langle #1 , #2 \right\rangle}

\newcommand{\R}[0]{ \IR\cup \{\infty \} }

\newcommand{\co}[1]{ #1^{\prime}}

\newcommand{\p}[0]{ p^{\prime}}

\newcommand{\m}[1]{   \mathcal{ #1 }}

\newcommand{ \W}[0]{ \mathcal{W}}

\newcommand{ \A}[1]{ \left\| #1 \right\|_H }

\newcommand{\B}[2]{ \left( #1 , #2 \right)_H }

\newcommand{\C}[2]{ \left\langle #1 , #2 \right\rangle_{  H^* , H } }

 \newcommand{\HON}[1]{ \| #1 \|_{ H^1} }

\newcommand{ \Om }{ \Omega}

\newcommand{ \pOm}{\partial \Omega}

\newcommand{\D}{ \mathcal{D} \left( \Omega \right)}

\newcommand{\DP}{ \mathcal{D}^{\prime} \left( \Omega \right)  }

\newcommand{\DPP}[2]{   \left\langle #1 , #2 \right\rangle_{  \mathcal{D}^{\prime}, \mathcal{D} }}

\newcommand{\PHH}[2]{    \left\langle #1 , #2 \right\rangle_{    \left(H^1 \right)^*  ,  H^1   }    }

\newcommand{\PHO}[2]{  \left\langle #1 , #2 \right\rangle_{  H^{-1}  , H_0^1  }}

 \newcommand{\HO}{ H^1 \left( \Omega \right)}

\newcommand{\HOO}{ H_0^1 \left( \Omega \right) }

\newcommand{\CC}{C_c^\infty\left(\Omega \right) }

\newcommand{\N}[1]{ \left\| #1\right\|_{ H_0^1  }  }

\newcommand{\IN}[2]{ \left(#1,#2\right)_{  H_0^1} }

\newcommand{\INI}[2]{ \left( #1 ,#2 \right)_ { H^1}}

\newcommand{\HH}{   H^1 \left( \Omega \right)^* }

\newcommand{\HL}{ H^{-1} \left( \Omega \right) }

\newcommand{\HS}[1]{ \| #1 \|_{H^*}}

\newcommand{\HSI}[2]{ \left( #1 , #2 \right)_{ H^*}}

\newcommand{\WO}{ W_0^{1,p}}
\newcommand{\w}[1]{ \| #1 \|_{W_0^{1,p}}}

\newcommand{\ww}{(W_0^{1,p})^*}

\newcommand{\Ov}{ \overline{\Omega}}

\def\IR{{\mathbb{R}}}

\title{Regularity of the extremal solutions associated to  elliptic systems}
\maketitle


\vspace{3mm}

\begin{abstract}
We examine the elliptic system given by
\begin{eqnarray*}
\qquad  \left\{ \begin{array}{lcl}
   -\Delta u =\lambda f(v) \quad \mbox{ in }  \Omega  \\
 -\Delta v =\gamma f(u)   \quad  \mbox{ in } \Omega,  \\
 u=v =0,  \quad \mbox{ on }  \pOm
\end{array}\right.
  \end{eqnarray*}
where $\lambda,\gamma$ are positive parameters, $\Omega$ is a smooth bounded
domain in $\IR^N$  and $f$ is a $C^{2}$ positive, nondecreasing and
convex function in $[0,\infty)$ such that $\frac{f(t)}{t}\rightarrow\infty$ as $t\rightarrow\infty$. Assuming
 $$0<\tau_{-}:=\liminf_{t\rightarrow\infty} \frac{f(t)f''(t)}{f'(t)^{2}}\leq \tau_{+}:=\limsup_{t\rightarrow\infty} \frac{f(t)f''(t)}{f'(t)^{2}}\leq 2,$$
 we show that the extremal solution $(u^*, v^*)$ associated to the above system is smooth provided\\ $N<\frac{2\alpha_{*}(2-\tau_{+})+2\tau_{+}}{\tau_{+}}\max\{1,\tau_{+}\}$, where $\alpha_{*}>1$ denotes  the largest root of the $2^{nd}$ order polynomial
$$P_{f}(\alpha,\tau_{-},\tau_{+}):=(2-\tau_{-})^{2} \alpha^{2}- 4(2-\tau_{+})\alpha+4(1-\tau_{+}).$$
As  a consequences, $u^*, v^*\in L^\infty(\Omega)$ for $N<5$. Moreover, if $\tau_{-}=\tau_{+}$, then  $u^*, v^*\in L^\infty(\Omega)$ for $N<10$.
\end{abstract}


\section{Introduction}

In this short note we examine the boundedness of the extremal solutions to the following system of equations:

 \begin{eqnarray*}
(P)_{\lambda,\gamma}\qquad  \left\{ \begin{array}{lcl}
\hfill   -\Delta u    &=& \lambda f(v) \quad \mbox{ in } \Omega  \\
\hfill -\Delta v &=& \gamma f(u)   \quad \mbox{ in }  \Omega,  \\
\hfill u &=& v =0 \quad  \mbox{ on }  \pOm,
\end{array}\right.
  \end{eqnarray*}    where $\Omega$ is a bounded domain in $\IR^N$ and $ \lambda, \gamma >0$ are positive parameters.  The nonlinearity $f$  satisfies
 \[ (R) \;\;  \mbox{$f$ is smooth, increasing and convex with $ f(0)=1 $ and $ f$ superlinear at $ \infty$.}\]

 Define $\mathcal{Q}:=\{(\lambda,\gamma),~\lambda,\gamma>0\}$,
 $$\mathcal{U}:=\{(\lambda,\gamma)\in\mathcal{Q}:~\text{there exists a smooth solution}~(u,v)~\text{of}~(P)_{\lambda,\gamma}\},$$
and set $\Upsilon:=\partial\mathcal{U}\cap\mathcal{Q} $. M. Montenegro in \cite{Mon} ( for a  more general system than $(P)_{\lambda,\gamma}$) showed that $\mathcal{U}\neq\emptyset$ and for every $(\lambda,\gamma)\in\mathcal{U}$ the problem $(P)_{\lambda,\gamma}$ has a minimal solution. Then, using monotonicity, for each $(\lambda^*,\gamma^*)\in\Upsilon$ one can define the extremal solution $(u^*,v^*)$ as a pointwise limit of minimal solutions of $(P)_{\lambda,\sigma\lambda}$ with $\sigma:=\frac{\gamma^*}{\lambda^*}$, which is always a weak solution to $(P)_{\lambda^*,\gamma^*}$. Moreover, for a $(\lambda,\gamma)\in\mathcal{U}$, the minimal solution $(u, v) ~\text{of}~ (P)_{\lambda,\gamma}$ is semi-stable in the sense that
there are  constants $\eta>0,~\zeta\geq0$ and $\chi\in H^1_0(\Omega)$ such that
\begin{equation}
-\Delta\zeta=\lambda f'(v)\chi+\eta\zeta,~~-\Delta\chi=\gamma g'(v)\zeta+\eta\chi,~~\text{in}~\Omega.
\end{equation}
For the proof  see \cite{Mon} (see also \cite{Craig1} for an alternative proof).\\
 In \cite{Mon} it is left open the question of the regularity of extremal solution $ (u^*,v^*)$. In the case when $f(t)=e^t$, in \cite{Craig} Cowan proved the extremal solutions  to $(P)_{\lambda,\sigma\lambda}$ are smooth  for $1 \leq N \leq 9$ under the further assumption $\frac{N-2}{8}<\frac{\gamma}{\lambda}<\frac{8}{N-2}$, and
 Dupaigne,  Farina and  Sirakov in \cite{DFS} proved it  without this restriction. The same result is also obtained by D\'avila and Goubet \cite{DG}.
 Furthermore, they proved that for $N\geq 10$,  the
singular set of any extremal solution of the system $(P)_{\lambda,\gamma}$ has Hausdorff dimension at most $N-10$.
 We now mention that some of the motivation for our proof of Theorem 1 in the  current paper comes  from the work of Dupaigne,  Farina and  Sirakov \cite{DFS}.\\

Now define
\begin{equation}
\tau_{-}:=\liminf_{t\rightarrow\infty} \frac{f(t)f''(t)}{f'(t)^{2}}\leq \tau_{+}:=\limsup_{t\rightarrow\infty} \frac{f(t)f''(t)}{f'(t)^{2}}.
\end{equation}
Our main result is the following.
\begin{thm}
 Let $f$  satisfy $(R)$ with $0<\tau_{-}\leq\tau_{+}<2$, and $\Omega$ an arbitrary bounded smooth domain. Also, let $ (u^*,v^*)$ denote the extremal solution associated with $ (P)_{\lambda,\gamma}$. Then $u^*, v^*\in L^\infty(\Omega)$ for
\begin{equation}
n<N(f):=\frac{2\alpha_{*}(2-\tau_{+})+2\tau_{+}}{\tau_{+}}\max\{1,\tau_{+}\}
\end{equation}
where $\alpha_{*}>1$ denotes  the largest root of the $2^{nd}$ order polynomial
\begin{equation}
P_{f}(\alpha,\tau_{-},\tau_{+}):=(2-\tau_{-})^{2} \alpha^{2}- 4(2-\tau_{+})\alpha+4(1-\tau_{+}).
\end{equation}
As  consequences,\\
i)  $u^*, v^*\in L^\infty(\Omega)$ for $N<5$.\\
i) If $\tau_{-}=\tau_{+}:=\tau$, then  $u^*, v^*\in L^\infty(\Omega)$ for $N<10$. Indeed, in this case we have
$$N(f)=2+4\frac{1+\sqrt{\tau}}{\tau}\geq 10.$$
\end{thm}

For example consider problem  $ (P)_{\lambda,\gamma}$ with  $f(t)=e^{t}$ or $e^{t^{\alpha}}$ ($\alpha>0$), then $\tau_{+}=\tau_{-}=1$, hence by Theorem 1, $u^*,v^*\in L^{\infty}(\Omega)$ for $n< 10$. The same is true for $f(u)=(1+u)^{p}$ ($p>1$) as in this case  we have  $\tau_{+}=\tau_{-}=\frac{p-1}{p}$. More precisely in the later case we have $u^*,v^*\in L^{\infty}(\Omega)$ for
$$n<2+\frac{4}{p-1}(p+\sqrt{p^2-p}).$$
This is exactly the same as the result  obtained in \cite{Craig1} and \cite{Haj} (corresponds to $p=\theta$ according to their notation).

\section{Preliminary estimates}
To prove the main result we use the following semistability inequality. For the proof see  \cite{CG,DFS}.
\begin{lemma} Let $(u,v)$ denote a semi-stable solution of $(P)_{\lambda,\gamma}$.  Then
\begin{equation}
\sqrt{\lambda \gamma}\int \sqrt{f'(u)f'(v)}\phi^2  \le  \int | \nabla \phi|^2,
\end{equation}
for all $\phi\in H^{1}_{0}(\Omega)$.
\end{lemma}
We need also the following lemmas.
\begin{lemma} Assume $\lambda\geq \gamma$. Then for any smooth solution to the  system $P_{\lambda,\gamma}$ we have
$$v\leq u \leq \frac{\lambda}{\gamma}v.$$

\end{lemma}
\textbf{Proof.} Take $w=u-v$. Then $w=0$ on $\partial\Omega$ and
$$-\Delta w=\lambda f(v)-\gamma f(u)\geq\lambda f(v)-\lambda f(u)=-\lambda \frac{f(u)-f(v)}{u-v}w:=-\lambda a(x)w,$$
where $a(x)=\frac{f(u)-f(v)}{u-v}\geq0$ because $f$ is increasing. Then by the maximum principle $w\geq0$ in $\Omega$.
Now take $\tilde{w}=\frac{\lambda}{\gamma}v-u$. Then $\tilde{w}=0$ on $\partial\Omega$ and using the above that $u\geq v$ we have
$$-\Delta \tilde{w}=\lambda f(u)-\lambda f(v)\geq0,$$ hence
$\tilde{w}\geq0$ in $\Omega$.\\
For the proof the next lemma we use the following standard regularity result, for the proof see Theorem 3 of \cite{Ser} and Theorems 4.1 and 4.3 of \cite{Tru}.
\begin{thm}
Let $u\in H^{1}_{0}(\Omega)$ be a weak solution of
\begin{equation}
\left\{\begin{array}{ll} \Delta u+c(x)u=g(x)& {\rm }\ x\in \Omega,\\~~u=0& {\rm }\ x\in \partial \Omega,
\end{array}\right.
\end{equation}
with $c,g\in L^{p}(\Omega)$ for some $p\geq1$.\\
Then there exists a positive constant $ C$  independent of $u$ such that if $p>\frac{n}{2}$ then
$$||u||_{L^{\infty}(\Omega)}\leq C(|u||_{L^{1}(\Omega)}+|g||_{L^{p}(\Omega)}).$$\\

\end{thm}
\begin{lemma} Assume for every semi-stable solution $(u,v)$  of $(P)_{\lambda,\gamma}$ with $\lambda\geq\gamma$ we have
$$||v||_{L^1(\Omega)}\leq C~~\text{and}~~||f'(v)||_{L^p(\Omega)}\leq C,$$
for some $p>\frac{N}{2}$, where $C$ is a constant independent of $(u,v)$. Then  $u^*, v^*\in L^\infty(\Omega)$.
\end{lemma}
\textbf{Proof.} We rewrite the first equation in $(P)_{\lambda,\gamma}$ as
$$\Delta u+\lambda \frac{f(v)-f(0)}{u}u=-\lambda f(0).$$
Taking $c(x):=\lambda \frac{f(v(x))-f(0)}{u(x)}$ then using Lemma 2 and the convexity of $f$ we have
$$0\leq c(x)\leq \lambda \frac{f(v(x))-f(0)}{v(x)}\leq \lambda^* f'(v).$$
Thus by the assumption and Theorem 2 we get $u^*\in L^\infty(\Omega)$, and by Lemma 2 we also get $v^*\in L^\infty(\Omega)$.
\section{Proof of Theorem 1}

\textbf{Proof of Theorem 1.}  Fix an $\alpha>1$ such that $P_{f}(\alpha,\tau_{-},\tau_{+})<0$. Such an $\alpha$ exists since we have $P_{f}(1,\tau_{-},\tau_{+})=(2-\tau_{-})^{2}-4<0$ and $P_{f}(+\infty,\tau_{-},\tau_{+})=+\infty$ .
Hence we can take positive numbers $\tau_{1}\in(0,\tau_{-})$  and $\tau_{2}\in(\tau_{+},2)$ such that
\begin{equation}
P_{f}(\alpha,\tau_{1},\tau_{2})<0.
\end{equation}

Now let $\phi(x)=\frac{\tilde{f}(u)^{\alpha}}{f'(u)^{\frac{\alpha}{2}}}$ in the semistabilty  inequality (5). Note that here for simplicity, we assumed that $f'(t)>0$ for $t>0$, this does not cause any problem, as in what follows we need only the behavior of $f$ and $f'$ at infinity. Then we get

\begin{equation}
\sqrt{\lambda \gamma}\int f'(u)^{\frac{1}{2}-\alpha}f'(v)^{\frac{1}{2}}\tilde{f}(u)^{2\alpha} \le  \lambda \int \theta(u)f(v),
\end{equation}
where
$$\theta(t)=\alpha^{2}\int_{0}^{t}\tilde{f}(s)^{2\alpha-2}f'(s)^{2-\alpha}\Big(1-\frac{\tilde{f}(s)f''(s)}{2f'(s)^{2}}\Big)^{2}ds.$$
First  we give an upper bound for the function $\theta$. By the definitions of $\tau_{\pm}$ there exists a $T>0$ such that $\tau_{1}\leq\frac{\tilde{f}(t)f''(t)}{f'(t)^{2}}\leq \tau_{2}$ for $t> T$ that also gives
\begin{equation}
0<1-\frac{\tau_{2}}{2}\leq 1- \frac{\tilde{f}(t)f''(t)}{2f'(t)^{2}}\leq 1-\frac{\tau_{1}}{2},~~\text{for}~ t> T .
\end{equation}
Using (9) we get
\begin{equation}
\theta(t)\leq \theta (T)+\alpha^{2} (1-\frac{\tau_{1}}{2})^{2}  \int_{T}^{t}\tilde{f}(s)^{2\alpha-2}f'(s)^{2-\alpha}ds,~~\text{for}~ t> T .
\end{equation}
Take $h(t):=\tilde{f}(t)^{2\alpha-1}f'(t)^{1-\alpha}$, then we have
$$h'(t)=(2\alpha-1)\tilde{f}(t)^{2\alpha-2}f'(t)^{2-\alpha}\Big(1- \frac{\alpha-1}{2\alpha-1} \frac{\tilde{f}(s)f''(s)}{f'(s)^{2}}\Big)$$
$$\geq (2\alpha-1)(1-\frac{\alpha-1}{2\alpha-1}\tau_{2})\tilde{f}(t)^{2\alpha-2}f'(t)^{2-\alpha},~~\text{for}~ t> T.$$
Using the above inequality in (10) we obtain
\begin{equation}
\theta(t)\leq C+ A  \tilde{f}(t)^{2\alpha-1}f'(t)^{1-\alpha},~\text{where}~A:=\frac{\alpha^{2}}{(2\alpha-1)} \frac{(1-\frac{\tau_{1}}{2})^{2}}{(1-\frac{\alpha-1}{2\alpha-1}\tau_{2})}~~\text{and}~C:=\theta (T)-Ah(T).
\end{equation}
Note that in the above we also used that $1-\frac{\alpha-1}{2\alpha-1}\tau_{2}>0$ which holds since $\tau_{2}<2$.
Now, the fact that the inequality $\frac{\tilde{f}(t)f''(t)}{f'(t)^{2}}\leq \tau_{2}$ for $t> T$ is equivalent to $\frac{d}{dt}(\frac{f'(t)}{\tilde{f}(t)^{\tau_{2}}})\leq0$ for $t> T$  gives
\begin{equation}
f'(t)\leq C_{1}\tilde{f}(t)^{\tau_{2}}~~ \text{for}~ t> T.
\end{equation}
Using this we obtain, for $t> T$
$$\tilde{f}(t)^{2\alpha-1}f'(t)^{1-\alpha}\geq f'(t)^{\frac{2\alpha-1}{\tau_{2}}-(\alpha-1)}\rightarrow\infty,~\text{as}~t\rightarrow \infty.$$
Now take an $\epsilon>0$. From the inequality above and (11), there exists an $T_{\epsilon}>T$ such that
\begin{equation}
\theta(t)\leq (A+\epsilon)\tilde{f}(t)^{2\alpha-1}f'(t)^{1-\alpha}, ~~\text{for}~ t>T_{\epsilon}.
\end{equation}
Also, we can find an $T'_{\epsilon}>0$ such that
\begin{equation}
f(t)\leq(1+\epsilon)\tilde{f}(t),~~\text{for}~ t>T'_{\epsilon}.
\end{equation}
Without loss of generality assume $\lambda\geq \gamma$ then from Lemma 2, $v\leq u\leq \frac{\lambda}{\gamma}v$. Using this,
taking $T''_{\epsilon}:=\max\{T_{\epsilon},T'_{\epsilon}\}$ and plugging (14), (13) in (8) we arrive at
\begin{equation}
\sqrt{\lambda \gamma}\int f'(u)^{\frac{1}{2}-\alpha}f'(v)^{\frac{1}{2}}\tilde{f}(u)^{2\alpha} \le  \lambda \Big(C_\epsilon+ (A+\epsilon)(1+\epsilon)\int_{v\geq T''_{\epsilon}} \tilde{f}(u)^{2\alpha-1}f'(u)^{1-\alpha}\tilde{f}(v)\Big),
\end{equation}
where
$$C_{\epsilon}:=\int_{v<T''_{\epsilon}} \tilde{f}(u)^{2\alpha-1}f'(u)^{1-\alpha}\tilde{f}(v)dx$$
is bounded by a constant independent of $u,v$, since  by Lemma 2 we have
 $$\{(u,v),~v\leq T''_{\epsilon}\}\subseteq [0,T''_{\epsilon}] \times [0,\frac{\lambda}{\gamma}~T''_{\epsilon}].$$
Letting
$$I:=\int f'(u)^{\frac{1}{2}-\alpha}f'(v)^{\frac{1}{2}}\tilde{f}(u)^{2\alpha}, $$
and replacing the integral on the right-hand side of   inequality (15) with integral over the full region $\Omega$ we get
\begin{equation}
\sqrt{\frac{\gamma}{\lambda }}I \le  C_\epsilon+ (A+\epsilon)(1+\epsilon)\int_{\Omega} \tilde{f}(u)^{2\alpha-1}f'(u)^{1-\alpha}\tilde{f}(v),
\end{equation}
By symmetry, taking

$$J:=\int f'(v)^{\frac{1}{2}-\alpha}f'(u)^{\frac{1}{2}}\tilde{f}(v)^{2\alpha},$$
we also get

\begin{equation}
\sqrt{\frac{\lambda}{\gamma }}J \le  C'_\epsilon+ (A+\epsilon)(1+\epsilon)\int_{\Omega} \tilde{f}(v)^{2\alpha-1}f'(v)^{1-\alpha}\tilde{f}(u),
\end{equation}
where $C'_\epsilon$ is bounded by a constant independent of $u,v$.\\
Now we write
$$\tilde{f}(u)^{2\alpha-1}f'(u)^{1-\alpha}\tilde{f}(v)=\Big(f'(u)^{\frac{1}{2}-\alpha}f'(v)^{\frac{1}{2}}\tilde{f}(u)^{2\alpha}\Big)^{\frac{2\alpha-1}{2\alpha}}
\Big(f'(v)^{\frac{1}{2}-\alpha}f'(u)^{\frac{1}{2}}\tilde{f}(v)^{2\alpha}\Big)^{\frac{1}{2\alpha}}.$$
Then, by the H\"{o}lder inequality we obtain
$$\int_{\Omega} \tilde{f}(u)^{2\alpha-1}f'(u)^{1-\alpha}\tilde{f}(v)\leq I^{\frac{2\alpha-1}{2\alpha}}J^{\frac{1}{2\alpha}}.$$
Using this in (16) we get
\begin{equation}
\sqrt{\frac{\gamma}{\lambda }}I \le  C_\epsilon+ (A+\epsilon)(1+\epsilon)I^{\frac{2\alpha-1}{2\alpha}}J^{\frac{1}{2\alpha}},
\end{equation}
and similarly from (17)
\begin{equation}
\sqrt{\frac{\lambda}{\gamma }}J \le  C'_\epsilon+ (A+\epsilon)(1+\epsilon)J^{\frac{2\alpha-1}{2\alpha}}I^{\frac{1}{2\alpha}}.
\end{equation}
Multiplying  inequalities (18) and (19), we get
\begin{equation}
\Big(1-(A+\epsilon)(1+\epsilon)\Big)IJ\leq C''_{\epsilon} (1+I^{\frac{2\alpha-1}{2\alpha}}J^{\frac{1}{2\alpha}}+J^{\frac{2\alpha-1}{2\alpha}}I^{\frac{1}{2\alpha}})
\end{equation}
where $C''_\epsilon$ is bounded by a constant independent of $u,v$. From (20) we deduce that if both of $I$ and $J$ are unbounded then we must have
$(1-(A+\epsilon)(1+\epsilon)\geq0$ and since $\epsilon>0$ was arbitrary we get $A\leq1$, which is equivalent to $P_{f}(\alpha,\tau_{1},\tau_{2})\geq0$, a contradiction. Hence, we proved that
\begin{equation}
f'(u)^{\frac{1}{2}-\alpha}\tilde{f}(u)^{2\alpha}f'(v)^{\frac{1}{2}} \in L^{1}(\Omega)~~\text{or}~~f'(v)^{\frac{1}{2}-\alpha}\tilde{f}(v)^{2\alpha}f'(u)^{\frac{1}{2}}\in L^{1}(\Omega).
\end{equation}
with a uniform bound in $L^{1}(\Omega)$ independent of $u,v$.
Now, it is easy to see that by our choice of $\alpha$ and the assumption that $\tau_+<2$, the function $y(t):=f'(t)^{\frac{1}{2}-\alpha}\tilde{f}(t)^{2\alpha}$  is an increasing function for $t$ large. Indeed, we have
$$y'(t)=(\alpha-\frac{1}{2})\tilde{f}(t)^{2\alpha-1}f'(t)^{\frac{3}{2}-\alpha}\Big(\frac{4\alpha}{2\alpha-1}-\frac{f(t)f''(t)}{f'(t)^{2}}\Big)$$
$$\geq (\alpha-\frac{1}{2})\tilde{f}(t)^{2\alpha-1}f'(t)^{\frac{3}{2}-\alpha}~(2-\tau_+)>0,$$
for $t$ sufficiently large. Hence, from (21) and the fact that $u\geq v$ we get
\begin{equation}
f'(v)^{1-\alpha}\tilde{f}(v)^{2\alpha} \in L^{1}(\Omega).
\end{equation}
From the inequality (12) and $\alpha>1$ we get
$$f'(t)^{1-\alpha}\tilde{f}(t)^{2\alpha}\geq \tilde{f}(t)^{(2-\tau_2)\alpha+\tau_2},~~t>T,$$
and also
$$f'(t)^{1-\alpha}\tilde{f}(t)^{2\alpha}\geq f'(t)^{\frac{(2-\tau_2)\alpha+\tau_2}{\tau_2}} ,~~t>T.$$
Hence, from (22) together with the above two inequalities we deduce that $\tilde{f}(v)^{(2-\tau_2)\alpha+\tau_2}\in L^{1}(\Omega)$ and also
$f'(v)^{\frac{(2-\tau_2)\alpha+\tau_2}{\tau_2}}\in L^{1}(\Omega)$. Now by the help of lemma 3 and the standard
elliptic regularity we get  $u^*, v^*\in L^\infty(\Omega)$ for
\begin{equation}
N<\max\{2\alpha(2-\tau_{2})+2\tau_{2},\frac{2\alpha(2-\tau_{2})+2\tau_{2}}{\tau_2}\}=\frac{2\alpha(2-\tau_{2})+2\tau_{2}}{\tau_{2}}\max\{1,\tau_{2}\}.
\end{equation}
Since we can choose $\tau_{2}$ arbitrary close to $\tau_{+}$ and $\alpha$ near to the largest root of the polynomial $P_{f}$, then  (23) completes the proof of the first part.\\
To see the second part, first note that we always have (since $\alpha^*>1$)
$$N(f)>2\alpha_{*}(2-\tau_{+})+2\tau_{+}>2(2-\tau_{+})+2\tau_{+}=4.$$
Also, if $\tau_-=\tau_+:=\tau$ then
$$\alpha^*=\frac{2+2\sqrt{\tau}}{2-\tau}.$$
Hence, $N(f)=2+4\frac{1+\sqrt{\tau}}{\tau}.$ Thus, using the fact that $\tau\leq1$ (since we always have $\tau_-\leq1$) we get $N(f)\geq 10$.
\hfill $ \Box$

\end{document}